\documentclass[12pt]{amsart}
\usepackage{amsmath}
\usepackage{amsfonts}
\usepackage{amsthm}

\theoremstyle{plain} \numberwithin{equation}{section}
\newtheorem{Theorem}{Theorem}
\newtheorem{Lemma}[Theorem]{Lemma}
\newtheorem{Proposition}[Theorem]{Proposition}

\theoremstyle{remark}

\title[trigonometric polynomial  potentials]
{Convergence of spectral decompositions of Hill operators with
trigonometric polynomial  potentials}
\author{Plamen Djakov}
\author{Boris Mityagin}
\begin{document}
\address{Sabanci University, Orhanli,
34956 Tuzla, Istanbul, Turkey}
 \email{djakov@sabanciuniv.edu}
\address{Department of Mathematics,
The Ohio State University,
 231 West 18th Ave,
Columbus, OH 43210, USA} \email{mityagin.1@osu.edu}

\begin{abstract}
We consider the Hill operator
$$
Ly = - y^{\prime \prime} + v(x)y, \quad   0 \leq  x \leq \pi,
$$
subject to periodic or antiperiodic boundary conditions, with
potentials $v$ which are trigonometric polynomials with nonzero
coefficients, of the form

 (i)  $ ae^{-2ix} +be^{2ix}; $

(ii)  $ ae^{-2ix} +Be^{4ix};  $

(iii)  $ ae^{-2ix} +Ae^{-4ix}  + be^{2ix} +Be^{4ix}. $

Then the system of eigenfunctions and (at most finitely many)
associated functions is complete but it is not a basis in $L^2
([0,\pi], \mathbb{C})$ if  $|a| \neq |b| $  in the case (i),  if $|A|
\neq |B| $ and neither $-b^2/4B$ nor $-a^2/4A$ is an integer square
 in the case (iii), and it is never a basis  in the case
(ii) subject to periodic boundary conditions.
\end{abstract}

\thanks{B. Mityagin acknowledges the hospitality of Sabanci
University, May--June, 2009}

\maketitle



\section{Introduction}
Convergence of spectral decompositions of ordinary differential
operators with various boundary conditions {\em bc} is a classical
area of research and has a long history -- see the monographs
\cite{Ti58, Na69, LS, MRH}.

  In the present paper we consider the Hill operators $L=L_{bc}(v)$
with smooth $\pi$-periodic potentials $v$
    \begin{equation}
  \label{1.1}  Ly = - y^{\prime \prime} + v(x)y, \quad   0 \leq  x \leq \pi,
  \end{equation}
  subject to periodic ($Per^+$)
  or antiperiodic ($Per^-$) boundary conditions:
  $$
  Per^\pm : \quad y(\pi) = \pm y(0), \quad y^\prime (\pi) = \pm
y^\prime (0).
  $$
 See basics and details in \cite{MW69}.

 Of course, if $v $ is real-valued, then $L_{Per^\pm}(v)$
 is a self-adjoint operator with a discrete spectrum.
 The system of its eigenfunctions
   \begin{equation}
  \label{1.2}
  \Phi = \{\varphi_k: \;\;   L \varphi_k = \lambda_k  \varphi_k, \; \; \|\varphi_k\|=1\}
  \end{equation}
is orthonormal, and the spectral decompositions
  \begin{equation}
  \label{1.3}
  f = \sum_k  \langle  f,\varphi_k  \rangle \varphi_k
  \end{equation}
converge (unconditionally) in $L^2 ([0,\pi])$ for every $f \in L^2
([0,\pi]).$

 If $v$ is a complex-valued potential the picture becomes more
complicated. If the boundary conditions are strictly regular then the
system of eigenfunctions and associated functions (SEAF) is a Riesz
basis in $ L^2([0, \pi])$ as it has been shown in
\cite{Du58,DS71,Ke64,Mi62}; see more details and history  in
\cite{Min99,Min06}. However, $Per^+, Per^-$ are regular but not
strictly regular boundary conditions. In this case properly chosen
two-dimensional block-decompositions do converge as it has been shown
by A. Shkalikov \cite{Sh79,Sh82,Sh83} (even in a more general context
of ordinary differential operators of higher order). For certain
classes of potentials, there have been given sufficient and necessary
conditions on whether blocks could be split into (one-dimensional)
eigenfunction decompositions \cite{MK98,DV05,Ma06-2,VS09}. Maybe, in
2006 A. Makin \cite{Ma06-1} and the authors \cite[Thm 71]{DM15} gave
first examples of such potentials that SEAF for periodic or
antiperiodic boundary conditions is NOT a basis in $L^2([0, \pi])$
even though all but finitely many eigenvalues are simple. The
existence of such potentials indirectly follows from the recent
results in \cite{GT09} as well.

We will extend many constructions and results of $SEAF$ divergence to
1D Dirac operators in an oncoming  paper \cite{DM24}.

In this paper we analyze low degree trigonometric polynomials and
show that the spectral decompositions of $L_{Per^\pm}$  {\em diverge}
if we exclude some exceptional values of coefficients of these
polynomials.

For example, if
$$
v(x) = a e^{-2ix} + b e^{2ix}, \quad a, b \in \mathbb{C}\setminus
\{0\},
$$
the SEAF decompositions  converge if and only if $|a|=|b|.$

In Section 2 we give the necessary preliminaries and prove a general
criterion (in terms of the Fourier coefficients of the potential $v,$
see Theorem 1) which says whether the SEAF is (or is not) a basis in
$ L^2 ([0,\pi]).$ Our constructions from \cite{DM10} are used in an
essential way when analyzing SEAF  related to trigonometric
potentials in Sections~3--5.

\section{Preliminary results}

It is well known that the spectra of the operators $L_{Per^\pm}$ are
discrete, and the following localization formulas hold (see, for
example, \cite[Prop~1]{DM10}):
\begin{equation}
  \label{2.2}
 Sp \, (L_{Per^\pm} ) \subset \Pi_N \cup \bigcup_{n>N, \,
 n\in\Gamma_{\pm}} D_n,    \quad \# \{Sp\,(L_{Per^\pm} )\cap D_n \} =2,
\end{equation}
where $D_n = n^2 + D,\;  D = \{z: \, |z| < 1\}, \; \Gamma_{+}=
2\mathbb{N}, \;  \Gamma_{-} =2\mathbb{N} -1,$
\begin{equation}
  \label{2.3}
\Pi_N = \{z=x+iy \in \mathbb{C}: \; |x| < (N+1/2)^2, \; |y| < N \},
\quad N=N(v).
\end{equation}
In either case the spectral block decompositions
\begin{equation}
  \label{2.7}
g = S_N g +\sum_{n>N, \,  n\in\Gamma_{\pm}}  P_n  g ,\quad \forall
\,g \in L^2 ([0,\pi]),
\end{equation}
where
\begin{equation}
  \label{2.8}
S_N = \frac{1}{2 \pi i} \int_{\partial \Pi_N} (z-L_{Per^\pm})^{-1}
dz, \quad P_n = \frac{1}{2 \pi i} \int_{|z-n^2|=1}
(z-L_{Per^\pm})^{-1} dz,
\end{equation}
converge unconditionally in $L^2 ([0,\pi]).$     This is true even
if the $\pi$-periodic potential $v$ is singular, i.e., $v \in
H^{-1}_{loc} (\mathbb{R}),$  as A. Savchuk and A. Shkalikov showed
\cite{SS03}. An alternative proof is given in \cite{DM19}.

We are going to provide in Theorem \ref{thm0} below sufficient
conditions which guarantee for large enough $n$ that each disc $D_n $
contains exactly two simple eigenvalues, and a criterion when the
two-dimensional spectral blocks in (\ref{2.7}) could be split into
one-dimensional spectral blocks so that to get an unconditional basis
in $L^2 ([0,\pi]).$

We shall use the following notations (compare with \cite{DM10}).
  For each $n\in \mathbb{N}$  a {\em walk} $x$
  from $-n$ to $n$ or from $n$ to
$-n$ is defined through its {\em sequence of steps}
\begin{equation}
  \label{2.9}
x=(x(t))_{t=1}^{\nu+1}, \quad 1\leq \nu=\nu(x)<\infty,
\end{equation}
where, respectively,
\begin{equation}
  \label{2.10}
\sum_{t=1}^{\nu+1} x(t) = 2n   \quad \text{or}  \quad
\sum_{t=1}^{\nu+1} x(t) = -2n.
\end{equation}
A walk $x$ is called {\em admissible} if its {\em  vertices} $j(t) =
j(t,x)$ given, respectively,  by
\begin{equation}
  \label{2.11a}
j(0) = -n  \quad \text{or} \;\;j(0) = +n
\end{equation}
and
\begin{equation}
  \label{2.11}
j(t) =-n + \sum_{t=1}^t x(i) \quad   \text{or}\quad j(t) =
 n + \sum_{t=1}^t x(i), \quad 1\leq t \leq \nu+1,
\end{equation}
satisfy
\begin{equation}
  \label{2.12}
 j(t) \neq \pm n \quad \text{for} \;\; 1\leq t \leq \nu.
\end{equation}

Let $X_n$ and $Y_n$ be, respectively, the set of all admissible walks
from $-n$ to $n$ and the set of all admissible walks from $n$ to
$-n.$ For each walk $x\in X_n$ or $x\in Y_n $ we set
\begin{equation}
  \label{2.13}
h(x;z) = \frac{\prod_{t=1}^{k+1} V(x(t)) }{\prod_{t=1}^k [n^2 -
j(t)^2 +z]}
\end{equation}
where $ V(m), \; m\in 2\mathbb{Z}$ are the Fourier coefficients of
the potential $v(x)$ with respect to the system $e^{imx}, \; m\in
2\mathbb{Z}.$ We set also
\begin{equation}
  \label{2.14}
B^+ (n,z) = \sum_{x\in X_n}  h(x,z), \quad B^- (n,z) = \sum_{x\in
Y_n} h(x,z).
\end{equation}

\begin{Theorem}
\label{thm0}
 Suppose  $v \in L^2 ([0,\pi]). $ If
\begin{equation}
\label{a1} B^+ (n,0) \neq 0, \quad B^- (n,0)\neq 0
\end{equation}
and
\begin{equation}
\label{a2} \exists c>0 : \quad  c^{-1}|B^\pm (n,0)| \leq  |B^\pm
(n,z)| \leq c \, |B^\pm (n,0)| \quad \forall \, z\in D
\end{equation}
for all sufficiently large even $n$ (if $bc= Per^+$) or odd $n$ (if
$bc= Per^-$), then

(a) there is $N=N(v) $ such that for $n>N$ the operator
$L_{Per^\pm}(v)$ has exactly two simple periodic (for even $n$) or
antiperiodic (for odd $n$) eigenvalues in the disc $D_n =  n^2 + D;$

(b)  a system of normalized eigenfunctions and associated functions
of $L_{Per^\pm} (v) $ is a Riesz basis in $L^2 ([0,\pi])$ if and only
if
\begin{equation}
\label{a3} 0 < \alpha := \inf_{n>N} \frac{|B^- (n,0)|}{|B^+ (n,0)|}
\quad \text{and}  \quad \beta :=\sup_{n>N} \frac{|B^- (n,0)|}{|B^+
(n,0)|} < \infty,
\end{equation}
where  we take $\inf $ and $ \sup $ over even $n$ if $bc= Per^+$ and
over odd $n$ if $ bc= Per^-.$
\end{Theorem}

{\em Remarks.} 1) Notice, that by (a)  the SEAF of $L_{Per^\pm} (v) $
has at most finitely many associated functions.

2) To avoid any confusion, let us emphasize that in
Theorem~\ref{thm0} are stacked together two {\em independent}
theorems: one for the case of periodic boundary conditions $Per^+$
(where we consider only even $n$), and another one for the case of
antiperiodic  boundary conditions $Per^-$ (where we consider only odd
$n$).

\begin{proof}
By the spectra localization formulas (\ref{2.2})  the operator
$L_{Per^\pm}(w)$ has, for each $n>N,$ two periodic (for even $n$) or
antiperiodic (for odd $n$) eigenvalues in the disc $n^2 +D $ (counted
with multiplicity). Moreover, by \cite{DM15} (see Lemma 21 and
Section 2.2, in particular, formula (2.23) and the three lines which
follow), the number $\lambda = n^2 + z, \; z\in D, $ is an eigenvalue
of $L_{Per^\pm}(v)$ if and only if  $z $ satisfies the basic equation
\begin{equation}
\label{a5} (z-a(n,z;v))^2 = B^+ (n,z;v)B^- (n,z,v), \quad z\in D,
\end{equation}
where  $a(n,z;v), B^\pm (n,z;v)$ are analytic functions of $z$ and
$v$ defined for $|Re \, z|<n. $  Next we show that for large enough
$n$ the equation (\ref{a5}) has exactly two roots in $D$ if counted
with multiplicity.

In view of \cite[Prop 28]{DM15} we have
\begin{equation}
\label{a8} |a(n,z,v)| \leq \frac{C}{n}, \quad |B^\pm (n,z;v)- V(\pm
2n)| \leq \frac{C}{n} \quad \text{if} \; z\in D,
\end{equation}
where $C=C(\|v\|)$ and $V(\pm 2n) $ are the $\pm 2n $-th Fourier
coefficients of the potential $v.$

Consider the family of potentials $w_t = t\cdot v, \; t \in [0,1].$
Since $ V(\pm 2n) \to 0$ as $n \to \infty, $ the inequalities
(\ref{a8}) imply
\begin{equation}
\label{a9} \sup_{\|z\|\leq 1}|a(n,z,tv)| \to 0, \quad
\sup_{\|z\|\leq 1} |B^\pm (n,z;tv)| \to 0  \quad \text{as} \;  n\to
\infty
\end{equation}
uniformly for $t \in [0,1].  $

Consider the function  $$F_n (z,t) = (z-a(n,z;tv))^2 - B^+
(n,z;tv)B^- (n,z,tv), \quad t \in [0,1].$$ In view of (\ref{a9}),
for large enough $n,$ the function $F_n (z,t) $ does not vanish on
the unit circle $\partial D. $ Therefore, the number of zeroes of
the equation (\ref{a5}) considered with $w=tv $   is given by
$$
\mathcal{N} (t) =\frac{1}{2\pi i } \int_{\partial D}
\frac{F_n^\prime (\zeta,t)}{F_n(\zeta,t)}d\zeta.
$$
Since the function $\mathcal{N} (t), \; t \in [0,1],$ is continuous
and takes integer values, it is a constant, so we have $ \mathcal{N}
(1) = \mathcal{N} (0). $ On the other hand, for zero potential the
basic equation is reduced to $ z^2 = 0, $ i.e., $\mathcal{N} (0)=2.$
Thus, for sufficiently large $n, $  say $n>N_1 $ the equation
(\ref{a5}) has exactly two roots in $D,$  counted with
multiplicities.

So, we have proved  for $n>N_1 $  that $\lambda = n^2 + z, \; z \in
D,$ is a periodic or antiperiodic value of algebraic multiplicity 2
if and only if $z$ is a double root of (\ref{a5}). Thus, the number
$\lambda = n^2 + z, \; z \in D,$ is a periodic or antiperiodic value
of algebraic multiplicity 2 if and only if $z$ satisfies the system
of the equation (\ref{a5}) and
\begin{equation}
\label{a6} 2(z-a(n,z)) \left (1-\frac{d}{dz}a(n,z)\right )
=\frac{d}{dz}\left ( B^+ (n,z)B^- (n,z) \right ).
\end{equation}
Therefore,  Part (a) of the theorem will be proved if we show that
there are at most finitely many $n$ such that the system (\ref{a5}),
(\ref{a6}) has a solution $z\in D.$

If $z(n) \in D $ is a root of (\ref{a5}), then by (\ref{a9})
\begin{equation}
\label{a10a} |z(n)| \leq |a(n,z(n))| + |B^+ (n,z(n))B^-
(n,z(n))|^{1/2} \to 0.
\end{equation}
Therefore, there is $\tilde{N}_1 =\tilde{N}_1 (v)>N_1 $  such that
\begin{equation}
\label{a10} |z(n)| \leq 1/2, \quad  n>\tilde{N}_1.
\end{equation}

Suppose that $n>\tilde{N}_1 $ and $z_n^* \in D $ satisfies the system
(\ref{a5}), (\ref{a6}). By the Cauchy inequality for the first
derivative, the first inequality in (\ref{a8}) implies
\begin{equation}
\label{a11} \left | \frac{d}{dz}a(n,z_n^*)\right | \leq 2C/n,
\end{equation}
while (\ref{a2}) and (\ref{a10}) yield
\begin{equation}
\label{a12} \left | \frac{d}{dz}B^+ (n,z_n^*)B^- (n,z_n^*)\right |
\leq 2 \sup_{|z|\leq 1}|B^+ (n,z)B^- (n,z)|
\end{equation}
$$
\leq 2c^2 \left | B^+ (n,0)B^- (n,0)\right | \leq 2c^4 |B^+
(n,z_n^*)B^- (n,z_n^*| .
$$
By (\ref{a5}), we have $ |z_n^* -a(n,z_n^*|=|B^+ (n,z_n^*) B^-
(n,z_n^*)|^{1/2}. $ Therefore, by (\ref{a11}) and (\ref{a12}), the
equation (\ref{a6}) implies
$$
2|B^+ (n,z_n^*) B^- (n,z_n^*)|^{1/2}(1-2C/n)\leq 2c^4 \left | B^+
(n,z_n^*)B^- (n,z_n^*)\right |,
$$
so it follows that $$ 1-2C/n\leq c^4\left | B^+ (n,z_n^*)B^-
(n,z_n^*)\right |^{1/2}.$$ By (\ref{a8}), the right-hand side of the
latter inequality tends to zero, so that inequality fails for large
enough $n.$ Hence, increasing if necessary $N_2, $ we obtain for $n>
N_2$ that the operator $L_{Per^\pm} (v) $ has no double periodic or
antiperiodic eigenvalues, i.e., (a) holds. \vspace{3mm}

Next we prove part (b)  of the theorem. In view of
(\ref{a1})--(\ref{a3}), for large enough $n$ the analytic functions
$B^+ (n,z) $ and $ B^- (n,z)$ do not vanish if $z \in D. $
Therefore, there are appropriate branches of $\log z $ (which depend
on $n$ and the choice of $\pm$) defined on a neighborhood of $ B^\pm
(n, \overline{D}). $ We set
$$
\text{Log}^\pm \, \left (B^\pm (n,z) \right ) = \log |B^\pm (n,z)| +
i \varphi_n^\pm (z);
$$
then
\begin{equation}
\label{a14}
 B^\pm (n,z) = |B^\pm (n,z)| e^{i\varphi_n^\pm (z)} \quad \forall \,
z\in D,\;\;n\geq N_2 (v)
\end{equation}
and the square root $\sqrt{B^+ (n,z) B^- (n,z)}$ is well defined by
\begin{equation}
\label{a15} \sqrt{B^+ (n,z) B^- (n,z)}=
 |B^+ (n,z) B^- (n,z)|^{1/2}  e^{\frac{i}{2}[\varphi_n^+ (z)+ \varphi_n^-
 (z)]}.
\end{equation}

Let us mention that the functions $\varphi_n^\pm $ are uniformly
Lipschitz on $\frac{1}{2}D; $ more precisely
\begin{equation}
\label{a16} |\varphi_n^\pm (z_1) -\varphi_n^\pm (z_2) | \leq 2c^2
|z_1-z_2| \quad \text{for} \;\; z_1,z_2 \in \frac{1}{2}D.
\end{equation}
Indeed, from (\ref{a2}) and the Cauchy inequality for the first
derivative it follows, for $|z|< 1/2,$ that
$$  \left | \frac{d}{dz}  \text{Log}^{\pm} \,
\left (B^\pm (n,z) \right ) \right | = \frac{1}{|B^\pm (n,z)|} \left
| \frac{d}{dz} B^\pm (n,z)\right | $$
$$
\leq \frac{c}{|B^\pm (n,0)|} \cdot 2\sup_{|z| \leq 1}|B^\pm
(n,z)|\leq \frac{c}{|B^\pm (n,0)|} \cdot 2c \,|B^\pm (n,0)| =2c^2.
$$

Now the basic equation (\ref{a5}) splits into the following two
equations
\begin{eqnarray}
\label{a17} z=\zeta_n^+ (z):=  a(n,z) + \sqrt{B^+ (n,z) B^- (n,z)}, \\
\label{a18} z=\zeta_n^- (z):=  a(n,z) - \sqrt{B^+ (n,z) B^- (n,z)}.
\end{eqnarray}
 For large enough $n,$ each of the equations (\ref{a17})
and (\ref{a18}) has exactly one root in the disc $D.$ Indeed,  in
view of (\ref{a8}), the Cauchy inequality for the first derivative
implies
$$
\sup_{|z|\leq 1/2} \left | d \zeta_n^{\pm}/dz \right | \to 0 \quad
\text{as} \quad n \to \infty.
$$
Therefore, for large enough $n$ each of the functions $\zeta_n^\pm $
is a contraction on the disc $\frac{1}{2}D, $  which implies that
each of the equations (\ref{a17}) and (\ref{a18}) has at most one
root in the disc $\frac{1}{2}D.$

On the other hand, by Part (a) and (\ref{a10}), for large enough $n$
the basic equation has two simple roots in $ \frac{1}{2}D $ and no
root on $D \setminus \frac{1}{2}D,$  which implies that each of the
equations (\ref{a17}) and (\ref{a18}) has exactly one root in the
disc $\frac{1}{2}D$ and no root on $D \setminus \frac{1}{2}D.$

For large enough $n,$ let $z_1 (n) $ (respectively $z_2 (n) $) be
the only root of the equation (\ref{a17}) (respectively (\ref{a18}))
in the unit disc $D.$  Let $f=f(n) $ and $g=g(n) $ be corresponding
unit eigenvectors of the operator $L= L_{Per^\pm},$ i.e.,
$\|f(n)\|=\|g(n)\|=1 $ and
$$
Lf(n) = (n^2 + z_1(n)) f(n), \quad Lg(n) = (n^2 +z_2(n))g(n).
$$

Let $P_n $ be the Riesz projections defined by (\ref{2.8}), and let
$P_n^0 $ be the Riesz projections associated with the free operator.
We have (e.g., see Proposition 11 in \cite{DM15})
\begin{equation}
\label{a21}     \dim P_n = \dim P_n^0 =2, \quad \|P_n - P_n^0\| \leq
C/n.
\end{equation}

Each of the projections $P_n, \; n>N,$ could be written as a sum of
one-dimensional projections on the subspaces generated by $f(n) $ and
$g(n)$ so that
$$P_n = P^1_n + P^2_n, \quad  P_n^1 P_n^2=P_n^2 P_n^1 =0.  $$
An elementary calculation shows that
$$
\|P_n^1\|= \|P_n^2\|= (1-|\langle f(n),g(n) \rangle \|^2)^{-1/2}.
$$
Therefore, the system of normalized eigenfunctions and associated
functions will be a Riesz basis if and only if
\begin{equation}
\label{a22} \limsup_{n\to \infty} |\langle f(n), g(n) \rangle | <1.
\end{equation}
We set
$$
f^0(n)= P_n^0 f(n), \quad  g^0(n)= P_n^0 g(n).
$$
From (\ref{a21}) it follows
$$
\|f(n) - f^0(n) \| = \|(P_n - P_n^0)f(n) \| \leq \|P_n - P_n^0\|
\leq C/n
$$
and $\|g(n) -g^0 (n)\| \leq C/n, \; |\langle f(n)-f^0(n), g(n)-g^0
(n) \rangle| \leq C/n^2.$  Since $ \|f (n)\|^2 =\|f^0 (n)\|^2+
\|f(n) - f^0(n) \|^2 $ and $ \langle f(n), g(n) \rangle = \langle
f^0(n), g^0(n) \rangle + \langle f(n)-f^0(n), g(n)-g^0 (n) \rangle,
$  we get
\begin{equation}
\label{a23} \|f^0(n)\|, \,\|g^0(n)\|\to 1,  \quad \limsup_{n\to
\infty} |\langle f(n), g(n) \rangle |=\limsup_{n\to \infty} |\langle
f^0(n), g^0(n) \rangle |.
\end{equation}

Then, by \cite[Lemma 21]{DM15} (see formula (2.4)),  $f^0 (n)$ is an
eigenvector of the matrix $\begin{pmatrix} a(n,z_1) &  B^+ (n,z_1)
\\B^- (n,z_1) & a(n,z_1)
\end{pmatrix}$ corresponding to its eigenvalue $z_1=z_1 (n), $ i.e.,
$$ \begin{pmatrix} a(n,z_1)-z_1 &  B^+ (n,z_1) \\B^- (n,z_1) &
a(n,z_1)-z_1  \end{pmatrix} f^0 (n) = 0. $$ Therefore, $f^0(n) $ is
proportional to the vector $\left (1, \frac{z_1 - a(n,z_1)}{B^+
(n,z_1)} \right )^T.$ Taking into account (\ref{a14}), (\ref{a15})
and (\ref{a17}) we obtain
\begin{equation}
\label{a25} f^0 (n) = \frac{\|f^0 (n)\|}{\sqrt{1+ \left | \frac{B^-
(n,z_1)}{B^+ (n,z_1)} \right |}} \begin{pmatrix} 1\\ \left |
\frac{B^- (n,z_1)}{B^+ (n,z_1)} \right
|^{1/2}e^{\frac{i}{2}[\varphi_n^- (z_1) - \varphi_n^+ (z_1)]}
\end{pmatrix}.
\end{equation}
In an analogous way, from (\ref{a14}), (\ref{a15}) and (\ref{a18}) it
follows
\begin{equation}
\label{a26} g^0 (n) = \frac{\|g^0 (n)\|}{\sqrt{1+ \left | \frac{B^-
(n,z_2)}{B^+ (n,z_2)} \right |}} \begin{pmatrix} 1\\ -\left |
\frac{B^- (n,z_2)}{B^+ (n,z_2)} \right
|^{1/2}e^{\frac{i}{2}[\varphi_n^- (z_2) - \varphi_n^+ (z_2)]}
\end{pmatrix}.
\end{equation}
Now,  (\ref{a25}) and (\ref{a26}) imply
$$
\langle f^0 (n), g^0 (n) \rangle  =\|f^0 (n)\|\|g^0 (n)\| \frac
{1-\sqrt{\left | \frac{B^- (n,z_1)}{B^+ (n,z_1)} \right |}\sqrt{
\left | \frac{B^- (n,z_2)}{B^+ (n,z_2)} \right |} \, e^{i\psi_n }}
{\sqrt{1+ \left | \frac{B^- (n,z_1)}{B^+ (n,z_1)} \right |}\sqrt{1+
\left | \frac{B^- (n,z_2)}{B^+ (n,z_2)} \right |}},
$$
where
$$\psi_n = \frac{1}{2}([\varphi_n^- (z_1(n)) -
\varphi_n^+ (z_1(n))]-[\varphi_n^- (z_2(n)) - \varphi_n^+
(z_2(n))]).
$$
In view of (\ref{a10a}) we have $z_1 (n) \to 0 $ and $z_2 (n) \to 0
$ as $ n\to \infty, $ so by (\ref{a16}) it follows
\begin{equation}
\label{a31} \psi_n \to 0 \quad \text{as} \;\; n \to \infty.
\end{equation}

We have
$$
|\langle f^0 (n), g^0 (n) \rangle |^2=\|f^0 (n)\|\|g^0 (n)\| \cdot
\Pi_n ,$$ where
$$
\Pi_n = \frac{1+\left | \frac{B^- (n,z_1)}{B^+ (n,z_1)} \right |
\left | \frac{B^- (n,z_2)}{B^+ (n,z_2)} \right |-2\sqrt{\left |
\frac{B^- (n,z_1)}{B^+ (n,z_1)} \right |}\sqrt{\left | \frac{B^-
(n,z_2)}{B^+ (n,z_2)} \right |} \cos \psi_n}{\left (1+ \left |
\frac{B^- (n,z_1)}{B^+ (n,z_1)} \right | \right ) \left (1+ \left |
\frac{B^- (n,z_2)}{B^+ (n,z_2)} \right | \right )}.
$$
If (\ref{a3}) fails, then there is a subsequence $n_k \to \infty $
such that
$$
\left | \frac{B^- (n_k,0)}{B^+ (n_k,0)} \right | \to 0 \quad
\text{or} \quad \left | \frac{B^- (n_k,0)}{B^+ (n_k,0)} \right | \to
\infty,
$$
which implies, in view of (\ref{a2}),  $\Pi_{n_k} \to 1.$ Therefore,
by (\ref{a23}),
$$ \limsup_{n\to
\infty} |\langle f(n), g(n) \rangle |= 1, $$ i.e., (\ref{a22}) fails,
so the system of normalized eigenfunctions and associated functions
is not a (Riesz) basis.

Suppose (\ref{a3}) holds. From  (\ref{a31}) it follows  $\cos \psi_n
>0 $ for large enough $n,$  so taking into account that $\|f^0 (n)\|,
\|g^0 (n)\|\leq 1, $ we obtain
$$ |\langle f^0 (n), g^0 (n) \rangle
|^2 \leq \Pi_n \leq
 \frac{1+\left | \frac{B^- (n,z_1)}{B^+ (n,z_1)} \right |
\left | \frac{B^- (n,z_2)}{B^+ (n,z_2)} \right |}{\left (1+ \left |
\frac{B^- (n,z_1)}{B^+ (n,z_1)} \right | \right ) \left (1+ \left |
\frac{B^- (n,z_2)}{B^+ (n,z_2)} \right | \right )} \leq \delta < 1,
$$
where
$$
\delta= \max  \left \{  \frac{1+xy}{(1+x)(1+y)}: \;
\frac{\alpha}{c^2} \leq x,y \leq c^2  \beta \right \}.
$$
Now (\ref{a23}) implies that (\ref{a22}) holds, hence the system of
normalized eigenfunctions and associated functions is a (Riesz)
basis in $L^2 ([0,\pi]).$ The proof is complete.
\end{proof}

In the next sections we consider the following
 three families of trigonometric polynomial potentials
\begin{equation}
  \label{2.34}
v(x) = a e^{-2ix} + b e^{2ix},
\end{equation}
\begin{equation}
  \label{2.35}
v(x) = a e^{-2ix} + B e^{4ix},
\end{equation}
\begin{equation}
  \label{2.36}
v(x) = a e^{-2ix} + A e^{-4ix}+ b e^{2ix} +B e^{2ix}
\end{equation}
and give conditions when the SEAF is a basis, in terms  of the
coefficients of these polynomials (see, respectively, Sections
3--5). \vspace{3mm}

 In all cases we consider in Sections 3--5   a special role is
played by {\em forward}  and {\em backward} walks.  We say that  $x $
is a  {\em forward}  (respectively, {\em backward}) walk if all steps
are positive, $x(t)>0  $ (respectively, negative,  $x(t)<0 $). Let
$X_n^+ $  and $Y_n^-$ be, respectively,  the set of all admissible
forward
 walks and the set of all admissible backward walks.

 \begin{Lemma}
 \label{lem2.6}
 If $\xi \in X_n^+ $  or $\xi \in Y_n^-, $  then for large enough $n$
 and $|z|\leq 1 $
 \begin{equation}
  \label{2.38}
 h(\xi, z) =  h(\xi, 0) (1+ \tau_n),   \qquad   |\tau_n| \leq  \frac{4 \log n}{n}.
  \end{equation}
  \end{Lemma}

\begin{proof}
By  (\ref{2.13}),
$$
\tau_n = \frac{h(\xi,z)}{h(\xi,0)} -1 =  \prod_{t=1}^\nu
\frac{n^2-j(t)^2}{n^2-j(t)^2 +z} -1 = e^{-w_n} -1,
$$
where $w_n = \sum_{t=1}^\nu \log \left (1+ \frac{z}{n^2 -j(t)^2}
\right ).$  Since $\xi \in X^+_n $  or $\xi \in Y^-_n,$ all vertices
$j(t)= j(t,\xi)$ are distinct, $ -n < j(t) < n  $ for $1 \leq  t \leq
\nu. $ Therefore, by the inequality $ |\log (1+\zeta)| \leq
\sum_{k=1}^\infty |\zeta|^k \leq 2 |\zeta|  $ for $|\zeta|\leq 1/2, $
for large enough $n$ it follows
$$
|w_n| \leq \sum_{t=1}^\nu \frac{2|z|}{n^2 -j(t)^2} \leq
\sum_{k=1}^{n-1} \frac{2}{n^2 - (-n+2k)^2 } = \frac{1}{n}
\sum_{k=1}^{n-1} \frac{1}{k} \leq \frac{2\log n}{n} \leq \frac{1}{2}.
$$
On the other hand, if $|w|\leq 1/2$ then $|e^{-w} -1| \leq
\sum_{k=1}^\infty |w|^k \leq 2|w|,$ which implies  (\ref{2.38}).
\end{proof}

\section{Potential $v= a e^{-2ix} + b e^{2ix}$}

We follow the notations and definitions of {\em walks, steps,
vertices} and functions $h, B^\pm $  given in
(\ref{2.9})--(\ref{2.14}).  The Fourier coefficients of the potential
$v=a e^{-2ix} + b e^{2ix}$ are
\begin{equation}
\label{3.1} V(-2)= a, \quad V(2) = b, \quad V(m) = 0 \;\; \text{for}
\; m \neq \pm 2.
\end{equation}

Let us focus on $B^+ (n,z).$ We say that a walk $x$ is $v$-{\em
admissible}, if $x$ is admissible and its steps are equal to $\pm 2.
$ If $x$ has $p$ steps equal to 2 and $q$ steps equal to $-2,$ then
\begin{equation}
\label{3.2} 2p -2q = 2n, \quad \text{so} \quad p = n+q,
\end{equation}
and \begin{equation} \label{3.3} p+q = \nu+1.
\end{equation}
We set
\begin{equation}
\label{3.4} X_n (q) = \{\text{$v$-admissible} \; x\in X_n \;\;
\text{with} \;\; q \; \text{steps} =-2 \}.
\end{equation}
Notice, that every $v$-admissible walk from $-n$ to $n$ has vertices
only between $-n$  and $n,$ and we have $x(1)=x(2) = 2.$ If
\begin{equation}
\label{3.5} i = \min \{t: \; x(t)\cdot x(t+1) <0 \},
\end{equation}
then
\begin{equation}
\label{3.6} x(t) = 2 \quad \text{if} \;\; 1 \leq t \leq i, \quad
x(i+1)= -2.
\end{equation}
We perform a "surgery" on $x$ by removing the steps $x(i)$ and $
x(i+1) $ and constructing a walk $\xi \in M^+ (q-1) $ such that
\begin{equation}
\label{3.7} \nu (\xi) = \nu  -2, \quad \nu =\nu (x),
\end{equation}
and
\begin{equation}
\label{3.8}
\xi (t) = \begin{cases}  x(t)  &    \text{for}  \;\; 1 \leq t \leq i-1\\
 x(t+2)  &        \text{for}  \;\; i \leq t \leq \nu-1.  \end{cases}
\end{equation}
Then
\begin{equation}
\label{3.9}
j(t, \xi) = \begin{cases}  j(t,x)  &    \text{for}  \;\; 1 \leq t \leq i-1\\
 j(t+2,x)  &        \text{for}  \;\; i \leq t \leq \nu -2.  \end{cases}
\end{equation}
Now we have
\begin{equation}
\label{3.10} h(x,z) =\frac{V(x(i))V(x(i+1))}{(n^2-j(i,x)^2
+z)(n^2+j(i+1,x)^2 + z)}\times h(\xi,z).
\end{equation}
With $c= |ab| $ the identity  (\ref{3.10})
 implies for $|z|\leq 1$
\begin{equation}
\label{3.11} \forall \, x \in X_n (q) \;\; \exists \xi \in X_n
(q-1):\quad
 |h(x,z)|   \leq \frac{c}{n^2} |h(\xi,z)|.
\end{equation}
Repeating the same procedure  $q$ times we come to  the inequality
\begin{equation}
\label{3.12} |h(x,z)|  \leq \left (\frac{c}{n^2} \right )^q
|h(\xi^*,z)|, \quad \exists \xi^* \in X_n (0), \quad |z|\leq 1.
\end{equation}
But $X_n (0)$ has only one element, and its only walk $\xi^* $ has
its steps, $n$ of them, equal to $2,$ so
\begin{equation}
\label{3.13} j(t,\xi^*) = -n + 2t, \quad  0 \leq t \leq n.
\end{equation}
We evaluate $h(\xi^*,0)$  and  estimate $h(\xi^*,z)$  below.

Let us notice that by (\ref{3.2})--(\ref{3.4})
\begin{equation}
\label{3.14} \# X_n (q) \leq \binom{p+q}{q}  = \binom{n+2q}{q} \leq
\begin{cases}
 2^{3q}  &  \text{if} \; q > n\\   \frac{1}{q!} (3n)^q  &  \text{if} \; q \leq n.      \end{cases}
\end{equation}
Therefore,
\begin{equation}
\label{3.15}  \sum_{q\geq 1} \sum_{x \in X_n (q)} |h(x,z)| \leq
\sigma_1 (n) \cdot  |h(\xi^*,z)|,
\end{equation}
where $$ \sigma_1 (n) = \sum_{q\geq 1} \binom{n+2q}{q} \left (
\frac{c}{n^2} \right )^q = \sum_{q=1}^n   +  \sum_{n+1}^\infty \leq
\sum_1^\infty \frac{1}{q!}   \left ( \frac{3nc}{n^2} \right )^q +
\sum_n^\infty  \left ( \frac{8c}{n^2} \right )^q $$ $$ \leq
\frac{3c}{n} e^{3c/n} +  \left ( \frac{8c}{n^2} \right )^n
\frac{1}{1- 8c/n^2} = O(1/n).$$ Thus, for $|z|\leq 1 $  we obtain
\begin{equation}
\label{3.17} B^+ (n,z) = h(\xi^*, z) + \sum_{q\geq 1} \sum_{x \in X_n
(q)} h(x,z)= h(\xi^*, z) (1+ O(1/n)).
\end{equation}

By Lemma~\ref{lem2.6}, we have $ h(\xi^*, z)= h(\xi^*, 0)(1+O(\log n
/n)),$ which leads to the following.
\begin{Lemma}
\label{lem3.5}
\begin{equation}
\label{3.18} B^+ (n,z) = h(\xi^*, 0) (1+O(\log n /n)), \quad \quad
|z| \leq 1.
\end{equation}
\end{Lemma}
The structure (\ref{3.13}) of $\xi^* $ makes possible to evaluate
$h(\xi^*,0)$ explicitly.
\begin{Lemma}
\label{lem3.6}
\begin{equation}
\label{3.19}
 h(\xi^*,0) = 4(b/4)^n [(n-1)!]^{-2}.
 \end{equation}
\end{Lemma}

\begin{proof} Indeed,
by (\ref{2.13}),  (\ref{3.1}) and  (\ref{3.13}), it follows that
$$ h(\xi^*,0)= b^n \prod_{t=1}^{n-1} [n^2- (-n+2t)^2]^{-1} = b^n
\left ( \prod_{t=1}^{n-1}  2t(2n-2t) \right )^{-1} = $$ $$  =b^n
4^{-n+1} \left (\prod_{t=1}^{n-1} t \right )^{-2} = 4(b/4)^n
 [(n-1)!]^{-2}.$$
\end{proof}

Lemmas \ref{lem3.5} and \ref{lem3.6} imply the following.
\begin{Proposition}
\label{prop3.7}
\begin{equation}
\label{3.22} B^+ (n,z) = 4(b/4)^n [(n-1)!]^{-2} (1+ O(\log n/n)),
\quad  |z| \leq 1.
\end{equation}
\end{Proposition}

To evaluate $B^- (n) $ we need to change forward walks to backward
walks, $b$ to $a,$  etc., which leads to the following
\begin{Proposition}
\label{prop3.7a}
\begin{equation}
\label{3.24} B^- (n,z) = 4(a/4)^n [(n-1)!]^{-2} (1+ O(\log n/n)),
\quad  |z| \leq 1.
\end{equation}
\end{Proposition}

The formulas (\ref{3.22}) and (\ref{3.24}) yield (\ref{a1}) and
(\ref{a2}), so Part (a) of Theorem ~\ref{thm0} implies that all but
finitely many of the eigenvalues of the operators $L_{Per^\pm} $ are
simple. Moreover, the following holds.

\begin{Theorem}
\label{thm3.8} Let $\{\varphi_k \} $ be a system of eigenfunctions
and associated functions of the operator
\begin{equation}
\label{3.25} -d^2/dx^2 + ae^{-2ix} + b e^{2ix}, \quad a,b \neq 0,
\end{equation}
subject to periodic $(Per^+)$ or antiperiodic $(Per^-)$ boundary
conditions. Then the spectral decomposition
\begin{equation}
\label{3.26} f= \sum c_k (f) \varphi_k
\end{equation}
converges (unconditionally) in $L^2 ([0,\pi])$ for each $f\in L^2
([0,\pi])$ if and only if
\begin{equation}
\label{3.27}  |a| = |b|.
\end{equation}
\end{Theorem}

\begin{proof}
If  $bc = Per^+$ we use the formulas (\ref{3.22}) and (\ref{3.24})
for even $n,$  while for antiperiodic boundary conditions $bc =
Per^-$ we use the same formulas with odd $n.$  By
Propositions~\ref{prop3.7} and \ref{prop3.7a},
\begin{equation}
\label{3.28} \frac{B^-(n,z)}{B^+ (n,z)} = \frac{a^n}{b^n} \left ( 1+O
\left (\frac{\log n}{n} \right ) \right ).
\end{equation}
If $|a|=|b| \neq 0,$ then  (\ref{a3}) holds, so by Theorem
\ref{thm0} the system $\{\varphi_k \} $ is an unconditional basis in
$L^2 ([0,\pi]).$

 If $|a| \neq |b|, $   then (\ref{a3}) fails, so Theorem
\ref{thm0} implies that the system $\{\varphi_k \} $ is not a basis
in $L^2 ([0,\pi]).$
 \end{proof}

For  Examples (\ref{2.35}) and (\ref{2.36}) we use the same general
scheme but technical details in estimations of $B^\pm (n,z) $ become
more complicated and interesting.

\section{Potential $v = a e^{-2ix} + B e^{4ix}.$}

In this case
\begin{equation}
\label{4.0} V(-2) = a, \quad V(4) = B; \quad V(j) = 0 \quad
\text{for} \;\; j \neq -2, 4.
\end{equation}
There is no symmetry in the structure of $v$-admissible forward and
backward walks (i.e., one cannot transform a forward part into a
backward  one or vice versa by replacing positive steps with the same
size negative steps). Therefore, we need to evaluate $B^- (n)$ and
$B^+ (n)$ separately.

Now we consider only periodic boundary conditions $Per^+,$ so $n$ is
even (see in Section~6 comments about the case $bc
=Per^-$).\vspace{3mm}

1.  First we  estimate $B^- (n,z), \; n=2m. $   We consider
$v$-admissible walks $x$ from $n$  to $-n;$  then
\begin{equation}
\label{4.2} x(t) = -2  \;\; \text{or}  \;\; 4, \quad   1 \leq t \leq
\nu +1, \quad \sum_1^{\nu +1} x(t)= -2n,
\end{equation}
where $\nu = \nu (x). $ If $p $  is the number of steps equal to 4
and $q$ is the number of steps equal to -2, then we have
\begin{equation}
\label{4.4} -2q + 4p = -2n= - 4m, \quad \text{so} \;\; q-2p = n=2m.
\end{equation}
Then $q $ should be even, say $q=2r,$ and we have
\begin{equation}
\label{4.5} r= p+m, \quad p+q = \nu +1.
\end{equation}

If $p=0$ then we have $q= n,$ and there is only one walk $\xi_* $
from $n $ to $-n$ with $n$ steps equal to $-2.$  We want to compare,
for any walk $x,$   $h(x,z)$ with $h(\xi_*,z). $  To this end we do a
surgery of $x$ by removing at once a triple of consecutive steps
$-2,-2, +4$  and get  a walk $\tilde{x}$ with $ p-1 $ steps equal to
4 and $q-2 $ steps equal to $-2.$ After that we estimate the ratio
$|h(x,z)|/|h(\tilde{x},z)|, $ and proceed further with another
surgery, and so on.

Let us denote by $Y_n (p) $  the set of all $v$-admissible walks from
$n$ to $-n$  having $p$ steps equal to 4. Suppose $x \in Y_n (p), \;
n>5. $ Then $x$ has a triple of consecutive steps $-2,-2, 4.$ Indeed,
one can easily see that
$$
x(1) = x(2) = x(3) = -2,
$$
because otherwise  $ n $ would be an intermediate vertex with
necessity, which is not possible for admissible walks by
(\ref{2.12}). Set
\begin{equation}
\label{4.11} i = \min\{ t: \; x(t+1) =4 \}
\end{equation}
and define $\tilde{x} \in Y_n (p-1) $ as
\begin{equation}
\label{4.12} \tilde{x} (t) = \begin{cases} x(t) = -2, &      1 \leq t
\leq i- 2\\   x(t+3),     &    i-1  \leq t \leq \nu -2.
\end{cases}
\end{equation}
Then
$$
h(x,z) = \frac{V(x(i-1))V(x(i)) V(x(i+1))}{\prod_{i-1}^{i+1} (n^2
-j(t,x)^2 +z)} \times h(\tilde{x},z),
$$
so with $c= |a^2 B| $ it follows
\begin{equation}
\label{4.15} \forall \, x\in Y_n (p) \;\; \exists \, \tilde{x} \in
Y_n (p-1): \quad |h(x,z) | \leq  \frac{c}{n^3}   |h(\tilde{x},z)|,
\quad |z|\leq 1.
\end{equation}

Repeating the same procedure $p$ times we obtain the inequality
\begin{equation}
\label{4.16}
 |h(x,z) | \leq  \left ( \frac{c}{n^{3}}  \right )^p  |h(\xi_*,z)| \quad \forall \, x\in Y_n (p),
\quad |z|\leq 1,
\end{equation}
where $\xi_*$  is the only walk of $Y_n (0).$  We have
 \begin{equation}
\label{4.18} j(t,\xi^*) = n-2t, \quad 0 \leq t \leq n.
\end{equation}

Let us notice that by (\ref{4.4})
\begin{equation}
\label{4.19} \# Y_n (p) \leq \binom{p+q}{p} = \binom{3p+2m}{p}.
\end{equation}
Therefore, by (\ref{4.16}) and (\ref{4.19}) it follows that
\begin{equation}
\label{4.20} \sum_{p\geq 1} \sum_{x\in Y_n (p)} |h(x;z)| \leq
\sigma_2 (n) \cdot  |h(\xi_*,z)|,
\end{equation}
where
\begin{equation}
\label{4.21}
 \sigma_2 (n) = \sum_{p\geq 1} \binom{3p+2m}{p}
 \left (  \frac{c}{n^3}  \right )^p = O(1/n^2).
\end{equation}
Indeed, since
$$
\binom{3p+2m}{p} \leq \begin{cases} 2^{5p}  &    \text{if}  \;\; p >
m,  \\   \frac{1}{p!} (5m)^p  &  \text{if}  \;\; p \leq m,
\end{cases}
$$
we have
$$
 \sigma_2 (n)  \leq \sum_{p=1}^m  + \sum_{p=m+1}^\infty  \leq
 \sum_{p=1}^\infty \frac{1}{p!}  \left ( \frac{5c}{8m^2} \right )^p +
 \sum_{p=m+1}^\infty  \left ( \frac{4c}{m^3} \right )^p
$$
$$
\leq \frac{5c}{8m^2} e^{\frac{5c}{8m^2}} + \left ( \frac{4c}{m^3}
\right )^{m+1} \frac{1}{ 1- \frac{4c}{m^3}} =O(1/n^2).
$$

Now we obtain, for $|z|\leq 1, $
\begin{equation}
\label{4.24} B^- (n,z) =h(\xi_*, z)+ \sum_{p\geq 1} \sum_{x\in Y_n
(p)} h(x;z)= h(\xi_*, z) (1+O(1/n^2)).
\end{equation}

By Lemma \ref{lem2.6}, $h(\xi_*, z)= h(\xi_*, 0)(1+ O(\log n/n)), $
which leads to the following.
\begin{Lemma}
\label{lem4.6}
\begin{equation}
\label{4.25} B^- (n,z) = h(\xi_*, 0) (1+ O(\log n/n)), \quad |z| \leq
1.
\end{equation}
\end{Lemma}
By (\ref{4.18}), we evaluate $h(\xi_*, 0) $ (compare with
Lemma~\ref{lem3.6}).
\begin{Lemma}
\label{lem4.7}
\begin{equation}
\label{4.26}
 h(\xi_*,0) = 4(a/4)^n [(n-1)!]^{-2}.
 \end{equation}
\end{Lemma}
In view of Lemmas  \ref{lem4.6} and \ref{lem4.7} the following holds.

\begin{Proposition}
\label{prop4.8}
\begin{equation}
\label{4.28} B^- (n,z) = 4(a/4)^n [(n-1)!]^{-2} \left (   1+ O(\log
n/n)    \right ), \quad |z| \leq 1.
 \end{equation}
\end{Proposition}
\vspace{2mm}

2.   To estimate $B^+ (n,z) $  we need to consider $v$-admissible
walks from $-n $ to $n.$ Let $X_n (q) $  be the set of all such walks
that have $q$ steps equal to -2. Notice, that $X_n (q) = \emptyset $
if $q $ is odd because  (compare with (\ref{4.2})) $\; - 2q +4p = 2n
= 4m $ so $ q+2m = 2p.$

For  even $q, $   say $q=2r,$ every $x \in X_n (q) $ has $p =
(2n+2q)/4 = m+r $ steps of length 4. The number of elements of $X_n
(q) $ could be estimated as
\begin{equation}
\label{4.33} \# X_n (q) \leq \binom{p+q}{q} = \binom{m+3r}{2r} \leq
\begin{cases} 2^{4r}   & \text{if} \;\; r>m,    \\ \frac{1}{(2r)!}
(4m)^{2r}  & \text{if} \;\; r \leq m.
\end{cases}
 \end{equation}

 \begin{Lemma}
 \label{lem4.10} For large enough $n,$
  \begin{equation}
\label{4.34} \forall \, x \in X_n (q) \; \exists  \tilde{x} \in X_n
(q-2): \quad |h(x,z)| \leq (c/n)^{5/2} \cdot |h(\tilde{x},z)|, \;\;
|z|\leq 1,
 \end{equation}
 where the constant $c>0 $ depends on $a$ and $B.$
 \end{Lemma}

\begin{proof}
Fix $x \in X_n (q);$  then one of the following two cases holds.

{\em Case 1.}  There are three consecutive steps $x(i_1), x(i_1+1),
x(i_1+2)$ with zero sum, i.e., $x(i_1)+ x(i_1+1)+ x(i_1+2)=0;$

{\em Case 2.} There is no triple of consecutive steps with zero sum,
i.e., two steps equal to $-2$ and one step equal to $+4$ (any order).

In Case 1 we set
\begin{equation}
\label{a.1} \tilde{x}(t) = \begin{cases} x(t),   &   1 \leq t\leq
i_1-1,\\ x(t+3), & i_1 \leq t \leq \nu -2.
\end{cases}
\end{equation}
Then each vertex of $\tilde{x}$ is a vertex of $x$ but $x$ has in
addition the vertices $j(i_1,x), j(i_1+1,x), j(i_1+2,x).$ Therefore,
it follows that
 \begin{equation}
\label{a.2} h(x,z)  = \frac{V(x(i_1))V(x(i_1+1))V(x(i_1+2))}
{\prod_{i_1}^{i_1+2} [n^2 - j(t,x)^2 +z]} \cdot h(\tilde{x},z),
\end{equation}
so
\begin{equation}
\label{a.3} |h(x,z)|  \leq \frac{|a^2 B|}{n^3}  \cdot
|h(\tilde{x},z)| \quad \text{if} \quad |z|\leq 1.
\end{equation}

In Case 2, set
\begin{equation}
\label{a.4} i_1 = \min \{t: \; x(t) =-2\}, \quad i_2=\min \{t > i_1:
\;  x(t)=-2\}
\end{equation}
and
\begin{equation}
\label{a.5} \tilde{x}(t) = \begin{cases} x(t)=4,   &   1 \leq t\leq
i_1-1,\\ x(t+2)= 4, &  i_1 \leq t \leq i_2 -3,\\ x(t+3), & i_2 -2
\leq t \leq \nu -2.
\end{cases}
\end{equation}
Notice that $i_2-i_1 \geq 3$ (otherwise we are in Case 1). Moreover,
from (\ref{a.4}) and (\ref{a.5}) it follows that $$ j(t,\tilde{x}) =
\begin{cases}  j(t,x) &   1\leq t \leq i_1 -1,\\
j(t+2, x) + 2  &    i_1 \leq t \leq i_2-3\\ j(t+2,x) &  i_2 -2 \leq t
\leq \nu -2.
\end{cases}
$$ Therefore,
 $ \frac{h(x,z)}{h(\tilde{x},z)} = P_1 (z)\cdot
P_2 (z), $ where $$ P_1 (z)=   \frac{a^2 B}{[n^2 - j(i_1,x)^2 +z][n^2
- j(i_1+1,x)^2 +z][n^2 - j(i_2-1,x)^2 +z]}
 $$ and $$ P_2 (z) = \prod_{t=i_1}^{i_2-3}
\frac{n^2 -(-n+4t)^2 +z}{n^2 -(-n+4t-2)^2 +z}  $$ Obviously, we have
$|P_1 (z) | \leq a^2 B/n^3.$ On the other hand,
$$P_2 (z)
\leq  P_2 (0) (1+ 4\log n/n) = 2\sqrt{n}(1+ 4\log n/n) $$ by
Lemma~\ref{lem2.6} and Lemma~\ref{lemP} below. Thus, Lemma
\ref{lem4.10} holds with $c= 2|a^2 B |. $
 \end{proof}

\begin{Lemma}
\label{lemP}
\begin{equation}
\label{a.11} \prod_{t=i}^j \frac{n^2 - (-n+4t)^2}{n^2 -(-n+4t-2)^2}
\leq 2\sqrt{n}, \quad 1 \leq i \leq j < n/2.
\end{equation}
\end{Lemma}

\begin{proof} Since
$$  \frac{n^2 - (-n+4t)^2}{n^2 -(-n+4t-2)^2} =
\frac{4t(2n-4t)}{(4t-2)(2n-4t+2)} \leq \frac{2t}{2t-1}, $$ the
product in (\ref{a.11}) does not exceed $\prod_1^n \frac{2t}{2t-1}.$
Since $ 2t/(2t-1) \leq \sqrt{t/(t-1)}$  for $t>1, $ we obtain
$$\prod_1^n \frac{2t}{2t-1} \leq 2
\sqrt{\frac{2}{1}}\sqrt{\frac{3}{2}}\cdots
\sqrt{\frac{n}{n-1}}=2\sqrt{n}.$$

\end{proof}

Now let us find the asymptotics of $B^+ (n,z).$  If we iterate
(\ref{4.34}) $r$ times then it follows
  \begin{equation}
\label{4.36} |h(x,z)| \leq (2c/n)^{5r/2} \cdot |h(\xi^*,z)|  \quad
\forall \, x \in X_n (q),  \;\; |z| \leq 1,
 \end{equation}
where $\xi^* $ is the only walk in $X_n (0);$ all its steps are equal
to 4, so
\begin{equation}
\label{4.37} j(t, \xi^*) = -n + 4t, \quad 0 \leq t \leq m.
 \end{equation}
By (\ref{4.33}) and (\ref{4.34}),
\begin{equation}
\label{4.38} \sum_{r\geq 1} \sum_{x\in X_n (2r)} |h(x,z) |  \leq
\sigma_3 (n)\cdot |h(\xi^*,z)|,
\end{equation}
where $$ \sigma_3 (n) \leq \sum_{r=1}^\infty \binom{m+3r}{2r} \left (
\frac{c}{n}  \right )^{\frac{5r}{2}} \leq \sum_{r=1}^m + \sum_{r>m}
$$
$$ \leq \sum_{r=1}^\infty \frac{1}{(2r)!} \left ( \frac{4n^2
c^{5/2}}{n^{5/2}} \right )^r +\sum_{r=m+1}^\infty 2^{4r}  \left (
\frac{c}{n}  \right )^{\frac{5r}{2}} $$      $$ \leq
\frac{2c^{5/2}}{\sqrt{n}} \exp (4c^{5/2}/\sqrt{n})  + O\left (
\frac{16c^{5/2}}{n^{5/2}} \right )= O(1/\sqrt{n}).
$$

Therefore, since $|B^+ (n,z) - h(\xi^*,z)| $ does not exceed the
left-hand side of (\ref{4.38}),  we obtain, in view of Lemma
\ref{lem2.6}, the following.

\begin{Lemma}
\label{lem4.12}
\begin{equation}
\label{4.41} B^+ (n,z) = h(\xi^*,0) (1+O(1/\sqrt{n})),  \quad |z|
\leq 1.
\end{equation}
\end{Lemma}

Next we evaluate $h(\xi^*,0).$
\begin{Lemma}
\label{lem4.13}
\begin{equation}
\label{4.42}
 h(\xi^*,0) = 16 (B/16)^{m} ((m-1)!)^{-2}.
 \end{equation}
\end{Lemma}

\begin{proof}
By (\ref{4.37}), $$  h(\xi^*,0)= B^m \cdot P^{-1},$$ where    $$ P
=\prod_1^{m-1} [n^2 - (-n+4t)^2] =\prod_1^{m-1} 4t (4m-4t)=16^{m-1}
((m-1)!)^2. $$ This proves (\ref{4.42}).
\end{proof}

Lemmas \ref{lem4.12} and \ref{lem4.13} imply the following (compare
with Proposition~\ref{prop4.8}).

\begin{Proposition}
\label{prop4.14} For even $n=2m$
\begin{equation}
\label{4.46}
 B^+ (n,z)= 16 (B/16)^{m} ((m-1)!)^{-2} (1+O(1/\sqrt{n})).
 \end{equation}
\end{Proposition}

11.  Now we apply Theorem \ref{thm0} and obtain the following.
\begin{Theorem}
\label{thm4.15} (a) If
\begin{equation}
\label{4.47} v= ae^{-2ix} +B e^{4ix}, \quad a,B \neq 0,
 \end{equation}
then all but finitely many of the eigenvalues of the operator
$L_{Per^\pm } $ are simple.

(b) If $(\psi_k)$  is a system of eigenfunctions and associated
 functions of the operator $L_{Per^+} (v),$
then this system is complete in $L^2 ([0,\pi])$ but it is not a basis
 in $L^2 ([0,\pi]).$
\end{Theorem}

\begin{proof}
In view of (\ref{4.28}) and (\ref{4.46}), the conditions (\ref{a1})
and (\ref{a2}) in Theorem~\ref{thm0} hold for even $n.$  Therefore,
by Part (a) of Theorem~\ref{thm0}, the operator $L_{Per^+} $ has at
most finitely many multiple eigenvalues.

Let  $\{\psi_k\}$ be a system of normalized eigenfunctions and
associated functions of the operator $L_{Per^+}. $ By (\ref{4.28})
and (\ref{4.46}),
$$
\lim_{n \;even}  \frac{B^-(n,0)}{B^+(n,0)} = 0,
$$
so the condition (\ref{a3}) fails. Thus, by Part (b) of
Theorem~\ref{thm0}, the system $\{\psi_k\}$ is not a basis in $L^2
([0,\pi]). $ This completes the proof.

 \end{proof}

\section{Potential $v= a e^{-2ix} + b e^{2ix} +A e^{-4ix} + B e^{4ix}$}

Now we analyze  trigonometric polynomials with four nonzero
coefficients of the form
\begin{equation}
\label{5.1} v= a e^{-2ix} + b e^{2ix} +A e^{-4ix} + B e^{4ix}.
\end{equation}
Since the set
\begin{equation}
\label{5.2} \{ k: \; V(k) \neq 0 \}   = \{ -2, -4, 2, 4 \}
\end{equation}
is symmetric,  it is enough to find  the asymptotics of $B^+ (n,z)$
in terms of  the coefficients $a, b, A,  B.$  Then we may obtain the
asymptotics of $B^- (n) $ just by exchanging the roles of $a, A $ and
$b, B.$

In our paper \cite{DM10},  we found the asymptotic behavior of the
spectral gaps of one-dimensional Schr\"odinger operator with a two
term potential $v= a \cos 2x + b \cos 4x, \; $  where $a$  and $b$
are real and nonzero. There, an essential part of the analysis is
related to the asymptotic behavior of the sums $ \sum_{x \in X_n}  h
(x,z), $
 so the techniques or even explicitly stated results from
\cite[Section 5]{DM10} give us tools to obtain asymptotics for $B^+
(n).$

Let $X_n  $  be the set of all walks $x$ from $-n$ to $n$ that are
$v$-admissible,  i.e., $ x(t) \in \{-2, -4, 2, 4\}, $ (\ref{2.12})
hold, and we have
\begin{equation}
\label{5.3} \sum_1^{\nu +1}  x(t)  = 2n,
\end{equation}
and let $X^+_n $ be the set of all $v$-admissible forward walks from
$-n$ to $n.$

In the case analyzed in Sections 3  (i.e., when $v= ae^{-2ix} + b
e^{2ix}) $
 there was only one forward walk. But now we have many such walks;
 more precisely, if A(n) is the number of solutions
 of (\ref{5.3}) with $x(t) =2 $ or $4, $ then $A(1) =1, \, A(2) = 2$
and  $A(n+1) = A(n) + A(n-1), $ so
 \begin{equation}
\label{5.6} \# X_n^+ = A(n) = \frac{1}{\sqrt{5}} \left ( \frac{1+
\sqrt{5}}{2} \right )^{n+1} - \frac{1}{\sqrt{5}} \left ( \frac{1-
\sqrt{5}}{2} \right )^{n+1}
\end{equation}
(Fibonacci numbers, see \cite[Sect. 4.1]{Cam}). \vspace{3mm}

2. For convenience, we change the parameters in (\ref{5.1}) by
setting
 \begin{equation}
\label{5.10} A= -\alpha^2, \quad a= -2\tau \alpha, \quad B = -
\beta^2, \quad b= -2\sigma \beta.
\end{equation}
In these notations the following statement (which is proven in
\cite{DM10}) holds.

\begin{Lemma}
\label{lem5.1} For even $n$
 \begin{equation}
\label{5.12} \sum_{\xi \in X^+_n} h(\xi,0) =
\frac{4(\beta/2)^n}{((n-1)!)^2} \prod_{i=1}^{n/2} [\sigma^2 -
(2i-1)^2],
\end{equation}
and for odd $n$
 \begin{equation}
\label{5.13} \sum_{\xi \in X^+_n} h(\xi,0) =
-\frac{4(\beta/2)^n}{((n-1)!)^2} \; \sigma \prod_{i=1}^{(n-1)/2}
[\sigma^2 - (2i)^2].
\end{equation}
\end{Lemma}

Product representations of $\cos t  $  and $\sin t $  show that
(\ref{5.12}) and (\ref{5.13}) could be rewritten as
 \begin{equation}
\label{5.14} \sum_{\xi \in X^+_n} h(\xi,0) =
\frac{4(i\beta/2)^n}{((n-2)!!)^2} \, \cos \left ( \frac{\pi
\sigma}{2}
 \right ) (1+O(1/n))
\end{equation}
for even $n,$ and
 \begin{equation}
\label{5.15} \sum_{\xi \in X^+_n} h(\xi,0) =
\frac{4i(i\beta/2)^n}{((n-2)!!)^2} \;\frac{2}{\pi} \sin \left (
\frac{\pi \sigma}{2}  \right ) (1+O(1/n))
\end{equation}
for odd $n.$

By Lemma \ref{lem2.6}, for every $\xi \in X^+_n $ we have
 \begin{equation}
\label{5.16}  | h(\xi, z)- h(\xi,0)  | \leq \frac{4\log n}{n}
\,|h(\xi,0)|,  \quad |z|\leq 1.
\end{equation}
These inequalities enable us to consider $z=0$ instead of $z $ in our
analysis of $B^+ (n,z).$ \vspace{3mm}

3. Now we are dealing with the difficulties brought by the huge size
of $X^+_n$  (see (\ref{5.6})).

For $\xi \in X_n^+, $  let $X_{n,\xi} $ denote the set of all walks
$x \in X_n\setminus X_n^+ $ such that each vertex $j(t, \xi)$ is a
vertex of $x$ also, i.e., $ j(s,\xi) = j(t_s, x)$ for some $t_s. $
Then we have
 \begin{equation}
\label{5.17} X_n\setminus X_n^+ = \bigcup_{\xi \in X_n^+} X_{n,\xi}.
\end{equation}
Indeed, for $x=(x(t))_{1}^{\nu+1} \in X_n $ define $t_0 =0$ and
 \begin{equation}
\label{5.18} t_{s+1} =\min \{  t>t_s \;: \;\; j(t,x) > j(t_s,x) \},
\quad 0 \leq s < \tilde{\nu},
\end{equation}
where
 \begin{equation}
\label{5.21} \tilde{\nu} = \min \{s: j(t_s, x) = n-4 \;\; \text{or}
\; n-2 \}.
\end{equation}
Define $\xi $ by the formula
 \begin{equation}
\label{5.22} \xi (s) = \begin{cases} j(t_s, x) - j(t_{s-1},x) &
\text{for} \;\; 1 \leq s \leq \tilde{\nu},\\ n- j(t_{\tilde{\nu}},x)
& \text{for} \;\; s= \tilde{\nu} +1.
\end{cases}
\end{equation}
Then $\xi \in X_n^+, $ and by the construction $x\in X_{n,\xi}.$

For $\xi \in X_n^+ $  and $m \in \mathbb{N}, $ let $X_{n,\xi,m}$ be
the set of walks $x\in X_{n,\xi} $ such that $x$ has $m$ more steps
than $\xi,$  i.e.,
 \begin{equation}
\label{5.24} X_{n,\xi,m} = \{x \in X_{n,\xi} \; :  \;\; \nu (x) - \nu
(\xi) =m \}.
\end{equation}
Then we have
 \begin{equation}
\label{5.25} X_{n,\xi} = \bigcup_{m=1}^\infty X_{n,\xi,m}.
\end{equation}

For $\xi \in X_n^+ $ and any $m$-tuple $I=(i_1, \ldots, i_m) $ of
integers  $i_\beta \in n + 2\mathbb{Z} \setminus \{\pm n \},$  let $
X_{n,\xi} (I)$ be the set of all walks $x$ with $\nu (\xi)+1 +m $
steps such that $I=(i_1, \ldots, i_m) $ and the sequence of the
vertices of $\xi$  are complementary subsequences of the sequence of
the vertices of $x.$ Then

\begin{Lemma} In the above notations, we have
\label{lem5.4}
\begin{equation}
\label{5.28}
 \# X_{n,\xi} (I) \leq 5^m \quad \forall \, I=(i_1, \ldots, i_m).
\end{equation}
\end{Lemma}
This is Lemma 12 in \cite{DM10}.

The following is an analogue of Lemma 13 in \cite{DM10}.
\begin{Lemma}
\label{lem5.5} There exists $n_1$ such that for $n\geq n_1 $
\begin{equation}
\label{5.29}   \sum_{x\in X_{n,\xi}} |h(x,z)| \leq \frac{K\log n}{n}
|h(\xi,z)| \quad \forall  \, \xi \in X_n^+, \;\; |z|\leq 1,
\end{equation}
where
\begin{equation}
\label{5.30} K= 40 C^2,\quad
 C=1+ \frac{\max (|a|, |b|, |A|, |B|)} {\min (|a|, |b|,
|A|, |B|)}.
\end{equation}
\end{Lemma}

\begin{proof}
In view of (\ref{5.25}), it is enough to show that
\begin{equation}
\label{5.32} \sum_{x\in X_{n,\xi,m}} |h(x,z)| \leq \left (\frac{K\log
n}{n} \right )^m |h(\xi,z)|,
\end{equation}
with $K$ and $C$ defined by  (\ref{5.30}). Indeed, if (\ref{5.32})
holds, then with $n_1 $ chosen so that $(K\log n)/n \leq 1/2 $ we
would have $$ \sum_{X_{n,\xi}} \frac{|h(x,z)|}{|h(\xi,z)|} \leq
\sum_{m=1}^\infty
 \left (\frac{K\log n}{n} \right )^m \leq \frac{K \log n}{n}, \quad
 n>n_1.
$$  which implies (\ref{5.29}).

To prove (\ref{5.32}), we use the inequality
\begin{equation}
\label{5.36}
 \sum_{x\in X_{n,\xi,m}}
|h(x,z)| \leq \sum_{I} \sum_{x\in X_{n,\xi} (I)} |h(x,z)|,
\end{equation}
where the first sum is taken over all $m$-tuples $I$ of integers
$i_\beta \in n+2\mathbb{Z},$  $i_\beta \neq \pm n.$ Fix such
$m$-tuple $I=(i_1,\ldots, i_m);$ then for every  $x \in X_\xi (I) $
$$ \frac{|h(x,z)|}{|h(\xi,z)|} = \frac{\prod_{t=1}^\nu
V(x_t)}{\prod_{\alpha=1}^{\tilde{\nu}} V(\xi_\alpha )} \times
\frac{1}{\prod_1^m (n^2 - i_\beta^2 +z)}. $$ We can split the first
factor $P$ as
\begin{equation}
\label{5.37}
 P = \prod_{\alpha=1}^{\tilde{\nu}} \left
(\frac{1}{V(\xi_\alpha)} \prod_{1+t_{\alpha -1}}^{t_\alpha} V(x(t))
\right ) \equiv \prod_{\alpha=1}^{\tilde{\nu}} r(\alpha).
\end{equation}

Let $d_\alpha= t_\alpha -t_{\alpha -1}; $ then $\sum_\alpha (d_\alpha
-1) = m. $  If $d_\alpha =1, $ then the ratio $r(\alpha) $ in
(\ref{5.37}) equals 1. Otherwise  $d_\alpha \geq 2,$ so, by the
inequality $d_\alpha \leq 2(d_\alpha -1), $ it follows that
$$|r(\alpha)| \leq C^{d_\alpha} \leq (C^2)^{d_\alpha -1}, \quad
\alpha=1, \ldots, \tilde{\nu}, $$ which implies $$ |P|
=\prod_{\alpha=1}^{\tilde{\nu}} r(\alpha) \leq \left ( C^2 \right
)^{\sum (d_\alpha -1)} = C^{2m}. $$ Therefore, taking into account
that $$ |n^2 -i^2 +z|^{-1} \leq 2|n^2 - i^2 | \quad \text{if} \quad
i\neq \pm n, \;\; |z|\leq 1 $$ we obtain
\begin{equation}
\label{5.46} \frac{|h(x,z)|}{|h(\xi,z)|}  \leq \frac{(2C^2)^m}{|n^2 -
i_1^2|\cdots |n^2 - i_m^2|}.
\end{equation}

Now by Lemma \ref{lem5.4}
\begin{equation}
\label{5.47}  \sum_{X_{n,\xi} (I)} \frac{|h(x,z)|}{|h(\xi,z)|}  \leq
\frac{(10C^2)^m}{|n^2 - i_1^2|\cdots |n^2 - i_m^2|},
\end{equation}
and by (\ref{5.36}) and the elementary inequality $$ \sum_{i \neq \pm
n} \frac{1}{|n^2-i^2|} \leq \frac{4 \log n}{n} \quad \text{for} \;\;
n\geq 10 $$ it follows that
\begin{equation}
\label{5.49}  \sum_{x\in X_{n,\xi,m}} \frac{|h(x,z)|}{|h(\xi,z)|}
\leq \sum_{i_1, \ldots, i_m \neq \pm n}\frac{(10C^2)^m}{|n^2 -
i_1^2|\cdots |n^2 - i_m^2|}
\end{equation}
$$ \leq (10C^2)^m  \left ( \frac{4 \log n}{n}  \right )^m = \left
( \frac{40 C^2 \log n}{n}  \right )^m. $$ Thus (\ref{5.32}) holds,
which completes the proof of Lemma~\ref{lem5.5}.
\end{proof}

4. Now we are going to complete the proof of the main result of this
section (compare with Step 5 and 6, pp. 187--190, in \cite{DM10}).

\begin{Proposition}
\label{prop5.6} If $\tau, \sigma $  given by (\ref{5.10}) are not
integers then for $|z| \leq 1$
\begin{equation}
\label{5.51} B^+ (n,z)= \frac{4(i\beta/2)^n}{((n-2)!!)^2} \, \cos
\left ( \frac{\pi \sigma}{2}
 \right ) \left( 1+O \left(\frac{\log n}{n}  \right) \right ),
\end{equation}
\begin{equation}
\label{5.51a} B^- (n,z)= \frac{4(i\alpha/2)^n}{((n-2)!!)^2} \, \cos
\left ( \frac{\pi \tau}{2}
 \right ) \left( 1+O \left(\frac{\log n}{n}  \right) \right ),
\end{equation}
for even $n,$ and
 \begin{equation}
\label{5.52} B^+ (n,z)= \frac{4i(i\beta/2)^n}{((n-2)!!)^2}
\;\frac{2}{\pi} \sin \left ( \frac{\pi \sigma}{2}  \right ) \left(
1+O \left(\frac{\log n}{n}  \right) \right ),
\end{equation}
 \begin{equation}
\label{5.52a} B^- (n,z)= \frac{4i(i\alpha/2)^n}{((n-2)!!)^2}
\;\frac{2}{\pi} \sin \left ( \frac{\pi \tau}{2}  \right ) \left( 1+O
\left(\frac{\log n}{n}  \right) \right ),
\end{equation}
for odd $n,$
with nonzero $\alpha, \beta, \tau, \sigma \in \mathbb{C} $  defined
in (\ref{5.10}).
\end{Proposition}

\begin{proof}
By symmetry of (\ref{5.2}), it is enough to prove only the estimates
for $B^+ (n,z).$

 From
(\ref{5.17}) and (\ref{5.29}) it follows that
 \begin{equation}
\label{5.55} \sum_{x\in X_n\setminus X_n^+} |h(x,z)| \leq \sum_{\xi
\in X_n^+} \sum_{x\in X_{n,\xi}}|h(x,z)| \leq K \frac{\log n}{n}
\sum_{\xi \in X_n^+} |h(\xi, z)|.
\end{equation}
Since $B^+ (n,z) = \sum_{x\in X_n} h(x,z) $ and $ |h(\xi, z)| \leq
 (1+ 4 \log n/n) |h(\xi,0)|$ due to Lemma~\ref{lem2.6}, the
 inequality (\ref{5.55}) implies
$$
 \left | B^+ (n,z) -\sum_{\xi \in X_n^+}h(x,z) \right |
\leq 2K \frac{\log n}{n} \sum_{\xi \in X_n^+} |h(\xi, 0)|.
$$
On the other hand, (\ref{5.16}) implies
$$\left |\sum_{\xi \in X_n^+}h(x,z) -\sum_{\xi \in X_n^+}h(x,0) \right |
\leq 4 \frac{\log n}{n} \sum_{\xi \in X_n^+} |h(\xi, 0)|.
$$
Therefore, we have
\begin{equation}
\label{5.57}  \left | B^+ (n,z) -\sum_{\xi \in X_n^+}h(x,0) \right |
\leq (2K+4) \frac{\log n}{n} \sum_{\xi \in X_n^+} |h(\xi, 0)|.
\end{equation}

Lemma \ref{lem5.1} gives an explicit formula for the sum $\sum_{\xi
\in X^+}h(\xi,0).$  The same formula could be used to find $\sum_{\xi
\in X^+}|h(\xi,0)|$ because $$ |h(\xi, 0)| = \frac{\prod_1^{\nu+1}
|V(\xi (t))|}{\prod_1^\nu (n^2 - j(t,\xi)^2)} = \frac{\prod_1^{\nu+1}
W(\xi (t))}{\prod_1^\nu (n^2 - j(t,\xi)^2)} \equiv h_w (\xi,0), $$
where $w$ is the potential $w =|a|e^{-2ix} +|b|e^{2ix} + |A| e^{-4ix}
+ |B| e^{4ix}$ with only nonzero Fourier coefficients
\begin{equation}
\label{5.61}
 W(-2)=|a|, \quad W(2)= |b|, \quad W(-4)=|A|, \quad
W(4)= |B|.
\end{equation}
Of course, now (\ref{5.10}) is replaced by
 \begin{equation}
\label{5.10a} |A|= -\tilde{\alpha}^2, \quad |a|= -2\tilde{\tau}
\tilde{\alpha}, \quad |B| = -  \tilde{\beta}^2, \quad b=
-2\tilde{\sigma} \tilde{\beta},
\end{equation}
with
\begin{equation}
\label{5.64} \tilde{\alpha}=i|\alpha|, \quad \tilde{\tau}=i|\tau |,
\quad \tilde{\beta}= i|\beta|, \quad \tilde{\sigma} =i|\sigma|,
\end{equation}
$\alpha, \beta, \tau, \sigma  $ coming from (\ref{5.10}).

Thus, we obtain,
\begin{equation}
\label{5.66} \sum_{\xi \in X_n^+} |h(\xi,0)| =
\frac{4(|\beta|/2)^n}{((n-2)!!)^2} \, \cosh \left ( \frac{\pi
|\sigma|}{2}
 \right ) \left( 1+O \left(\frac{1}{n}  \right) \right )
\end{equation}
for even $n,$ and
 \begin{equation}
\label{5.67} \sum_{\xi \in X_n^+} |h(\xi,0)|=
\frac{4(|\beta|/2)^n}{((n-2)!!)^2} \;\frac{2|\sigma|}{\pi} \sinh
\left ( \frac{\pi |\sigma|}{2}  \right ) \left( 1+O \left(\frac{1}{n}
\right) \right )
\end{equation}
for odd $n.$

Now we continue to analyze $B^+ (n,z).$ In view of (\ref{5.14}),
(\ref{5.15}), (\ref{5.66}) and  (\ref{5.67}) we have
\begin{equation}
\label{5.71}  \left |\sum_{\xi \in X_n^+} h(\xi,0)   \right |= \left
( \sum_{\xi \in X_n^+} |h(\xi,0)|\right )\cdot R^+ \cdot \left( 1+O
\left(\frac{1}{n} \right) \right ),
\end{equation}
where
\begin{equation}
\label{5.72} R^+ = \frac{|\cos\frac{\pi \sigma}{2}|}{\cosh\frac{\pi
|\sigma|}{2}} \quad \text{for even} \; n,  \quad  R^+ = \frac{|\sin
\frac{\pi \sigma}{2}|}{\sinh \frac{\pi |\sigma|}{2}} \quad \text{for
odd} \; n.
\end{equation}
Therefore, by (\ref{5.57}) and (\ref{5.71}) it follows that
\begin{equation}
\label{5.77}
 \left | B^+ (n,z) -\sum_{\xi \in X^+}h(x,0) \right
| \leq M \frac{\log n}{n} \left | \sum_{\xi \in X^+} h(\xi, 0)
\right|,
\end{equation}
where $M = (2K + 4)/R^+,$  so if $R^+ \neq 0,$ we have
\begin{equation}
\label{5.79}
 B^+ (n,z)= \left  ( \sum_{\xi \in X_n^+}
h(\xi, 0) \right ) \left (1+ O \left ( \frac{\log n}{n} \right )
\right ).
\end{equation}
The condition $R^+ \neq 0$ holds in both  $Per^\pm$ cases if and only
if $\sigma $ is not an integer. By (\ref{5.14}) and (\ref{5.15}), we
know the sum in the right-hand side of (\ref{5.79}). This completes
the proof of Proposition~\ref{prop5.6}.

\end{proof}

{\em Remark.} In analysis of $B^- (n,z) $ as an analog of  $R^+ $ we
would have
\begin{equation}
\label{5.72a} R^- = \frac{|\cos\frac{\pi \tau}{2}|}{\cosh\frac{\pi
|\tau |}{2}} \quad \text{for even} \; n, \quad  R^- = \frac{|\sin
\frac{\pi \tau}{2}|}{\sinh \frac{\pi |\tau|}{2}} \quad \text{for odd}
\; n.
\end{equation}
In terms of the coefficients $ a,b,A,B$ the condition  ``$\tau,
\sigma  $ are not integers`` in Proposition~\ref{prop5.6} holds if
and only if neither $-b^2/(4B),$  nor  $-a^2/(4A)$ is an integer
square. \vspace{3mm}

8. Now by the general scheme given in Theorem~\ref{thm0}, we obtain
the following.

\begin{Theorem}
\label{thm5.8} Consider the Hill operator $L_{Per^\pm} (v),$ where
\begin{equation}
\label{5.81} v=ae^{-2ix} + Ae^{-4ix} + be^{2ix}+Be^{4ix}
\end{equation}
with $a,b, A, B \neq 0$  and
\begin{equation}
\label{5.80}  \text{neither}\;\; -b^2/(4B), \;\; \text{nor} \;\;
-a^2/(4A) \quad \text{is an integer square}.
\end{equation}
All eigenvalues of $L_{Per^\pm} (v)$ but finitely many are simple;
the system $\Phi = \{\varphi_k \} $ of eigenfunctions and associated
functions is complete. $\Phi $ is an (unconditional) basis in $L^2
([0,\pi])$ if and only if $|A|=|B|.$
\end{Theorem}

\begin{proof}
In view of (\ref{5.80}), we may apply Proposition~\ref{prop5.6}.
Then, (\ref{5.51}) and (\ref{5.51a}) imply the conditions (\ref{a1})
and (\ref{a2}) in Theorem~\ref{thm0} for even $n,$ and (\ref{5.52})
and by (\ref{5.52a}) imply (\ref{a1}) and (\ref{a2}) for odd $n.$
Therefore, by Part (a) of Theorem~\ref{thm0}, each of the operators
$L_{Per^+} $ and $L_{Per^-} $ has at most finitely many multiple
eigenvalues.

Let  $\Phi=\{\varphi_k\}$ be a system of normalized eigenfunctions
and associated functions of the operator $L_{Per^+}. $ Then by
(\ref{5.51}) and (\ref{5.51a}) we have
\begin{equation}
\label{5.83} \frac{|B^- (n,0)|}{|B^+ (n,0)|} = \left |\frac{A}{B}
\right |^{n/2} \left |  \frac{\cos \frac{\pi \tau}{2}}{\cos \frac{\pi
\sigma}{2}} \right | \left ( 1+O \left( \frac{\log n}{n} \right )
\right ), \quad n\;\text{even}.
\end{equation}
Therefore, $$ \lim_{n \;even}  \frac{|B^-(n,0)|}{|B^+(n,0)|} =
\begin{cases}  0  &  |A|<|B|\\
\infty &   |A|>|B|\\
 \frac{\cos \frac{\pi \tau}{2}}{\cos \frac{\pi
\sigma}{2}}  &  |A|=|B| \end{cases}
$$
so the condition (\ref{a3}) fails if $|A|\neq |B|$ and holds if
$|A|=|B|.$ Thus, by Part (b) of Theorem~\ref{thm0}, if $|A|\neq |B|$
the system $\Phi_k$ is not a basis in $L^2 ([0,\pi]), $ and if $|A|=
|B|$ then $\Phi_k$ is an unconditional basis in $L^2 ([0,\pi]).$

In the same way, the conditions (\ref{5.52}) and (\ref{5.52a}) imply
the theorem for antiperiodic boundary conditions $Per^-.$ This
completes the proof.

\end{proof}

\section{Comments, conclusion}

1. In Section 4 we consider only periodic boundary conditions in the
case of potentials $v(x) = ae^{-2ix} + B e^{4ix}. $ In the case of
antiperiodic boundary conditions we need to analyze $B^\pm (n,z)$ for
odd $n.$ It turns out that most of the estimates done in Section 4
can be carried on for {\em odd} $n$ as well. But the crucial fact
\begin{equation}
\label{6.1} B^+ (n,z) = h(\xi^*,0) \left (1+O(\log n/n) \right ),
\quad n \;\; \text{even,}
\end{equation}
(see (\ref{4.41}), (\ref{4.42}), (\ref{4.46})) does not have a
reasonable analog if $n$ is odd. This observation and attempts to
follow the scheme which was successful for even $n$ are interesting
because they lead to some combinatorial problems and maybe give some
hints how the case $bc = Per^- $ could be studied.

Now, for an odd $n= 2m+1$ we write formulas that are analogous to
(\ref{4.2})--(\ref{4.5}).  Let $x = (x(t)_1^{\nu+1}$ be a
$v$-admissible walk from $-n$ to $n$ with $x(t) \in \{-2,4\}. $ We
denote by $p$ and $q,$ respectively,  the number of steps  equal to 4
and the number of steps equal to $-2.$  Then $ 4p-2q =2n =2(2m+1), $
so we have
$$ 2p = 2m + 1 + q, \quad p+ q = \nu + 1. $$ Now $q$ is
odd, say $q= 2r +1, $ and  $q=1$  is the minimal possible value of
$q.$

Let $X_n^+ (q) $  denote the set of all admissible walks with $q$
steps equal to $-2.$  By repeating the constructions of Section~4 one
may prove the following statements.

\begin{Lemma}
\label{lem6.1} If $r>0, $ then
\begin{equation}
\label{6.7} \sum_{x \in X_n^+ (2r+1)} |h(x,z)| \leq  \left (
\frac{c}{n^{5/2}} \right )^r \sum_{\xi \in X_n^+ (1)} |h(\xi,z)|.
\end{equation}
\end{Lemma}

\begin{Lemma}
\label{lem6.2} For large enough $n$
\begin{equation}
\label{6.8} \left | B^+ (n,z) - \sum_{\xi \in X_n^+ (1)} h(\xi,z)
\right | \leq \frac{c}{n^{5/2}} \sum_{\xi \in X_n^+ (1)} |h(\xi,z)|,
\quad |z|\leq 1.
\end{equation}
\end{Lemma}

\begin{Lemma}
\label{lem6.3} For large enough $n$
\begin{equation}
\label{6.9} \left | \sum_{\xi \in X_n^+ (1)} h(\xi,z) -\sum_{\xi \in
X_n^+ (1)} h(\xi,0)\right | \leq  4\frac{\log n}{n}  \sum_{\xi \in
X_n^+ (1)} |h(\xi,0)|
\end{equation}
\end{Lemma}

So far it is OK. But
\begin{equation}
\label{6.11} \# X_n^+ (1) = m+2,
\end{equation}
and in order to apply (\ref{6.8}) and (\ref{6.9}) we need to evaluate
$$ H_0^* =\sum_{\xi \in X_n^+ (1)} |h(\xi,0)| \quad \text{and} \quad H_0
= \sum_{\xi \in X_n^+ (1)} h(\xi,0)$$  and be sure that $H_0 \neq 0.
$

We can evaluate $H^* $ and $H_0 $ (see Proposition~\ref{prop6.4}) but
$H_0 = 0 $ -- see (\ref{6.24}).

Any walk $\xi \in X_n^+ (1) $ has only one step equal to $-2$ but
that step could appear on the left of $-n $ (denote that walk by
$\xi^- $ ), on the right of $n$ (denote that walk by $\xi^+ $) and
anywhere between $-n $ and $n$ (denote the set of all such walks by
$\tilde{X}_n^+(1)).$ With $p= m+1 $ the numerator in $h(\xi, 0 ) $ is
equal to $ab^{m+1}$ for every $\xi \in M^+ (1), $ so we can assume in
the calculations which follow that $a=b=1.$ Then the sum $H_0 $ has
two negative terms, namely
\begin{equation}
\label{6.14} h(\xi^-,0)= h(\xi^+ ,0) = 1/P,
\end{equation}
where, with $n=2m+1, $
\begin{equation}
\label{6.15} P = (n^2 - (-n-2)^2) \prod_{\tau=0}^{m-1} [n^2 -
(-n+2+4\tau )^2 ] = -2 (2n+2)\prod_{\tau=0}^{m-1} (2+4\tau)(4m-4\tau)
\end{equation}
$$ = - 8(m+1) 8^m m! (2m-1)!! = -8 (m+1) 4^m (2m)!$$ Therefore,
\begin{equation}
\label{6.16} h(\xi^-,0)+ h(\xi^+ ,0) = - \frac{1}{m+1} \cdot
\frac{1}{4^{m+1}} \cdot \frac{1}{(2m)!}
\end{equation}
The remaining walks $\xi \in \tilde{X}_n^+ $ give a sum of positive
terms of the form $(P(t) P(s))^{-1} $  with $s= m-t+1,$ where
\begin{equation}
\label{6.18} P(t)= \prod_{\tau=1}^t [n^2 - (-n+4\tau )^2 ]=
\prod_{\tau=1}^t 4\tau (4m+2 - 4\tau) =
\end{equation}
$$ 8^t t! \prod_{\tau=0}^t [2(m-\tau +1) -1] = 8^t t!
\frac{(2m-1)!!}{(2(m-t)-1)!!} \cdot \frac{2^{m-t} (m-t)!}{(2(m-t))!!}
\cdot \frac{(2m)!!}{2^m m!} $$ $$=4^t \, t! \frac{(2m)!}{m!}  \cdot
\frac{(m-t)!}{(2(m-t))!}.$$ Then with $t+s = m+1, \; 1 \leq t \leq m,
$ we have
\begin{equation}
\label{6.19} P(s)=4^{m+1-t} (m+1-t)! \frac{(t-1)!}{(2(t-1))!} \cdot
\cdot \frac{(2m)!}{m!}
\end{equation}
and
\begin{equation}
\label{6.20} P(t) P(s)=4^{m+1} \binom{2m}{m} (2m)! \cdot
\frac{t}{\binom{2(t-1)}{t-1}} \cdot
\frac{m+1-t}{\binom{2(m-t)}{m-t}}.
\end{equation}
Next, we use  Catalan numbers (see \cite[Section 4.5, (4.5.1) and
(4.5.2)]{Cam} or \cite[pp. 117, (14.10)--(14.12)]{LW})
\begin{equation}
\label{6.21} C_k = \frac{1}{k} \binom{2k-2}{k-1}, \quad k \geq 1,
\end{equation}
and the  fundamental recurrence for Catalan numbers
\begin{equation}
\label{6.23}
 C_{k+1} = \sum_{i=1}^{k} C_i C_{k+1-i}.
 \end{equation}
In view of (\ref{6.21}) and (\ref{6.23}), we obtain
\begin{equation}
\label{6.22} \sum_{t+1}^m (P(t) P(m+1-t))^{-1} = \frac{1}{4^{m+1}}
\cdot \frac{1}{m+1} \cdot \frac{1}{(2m)!} \left [C_{m+1}^{-1}
\sum_{t=1}^m C_{t}C_{m+1-t} \right ]
\end{equation}
$$=\frac{1}{4^{m+1}} \cdot \frac{1}{m+1} \cdot \frac{1}{(2m)!}. $$
Now (\ref{6.16}) and (\ref{6.22}) imply the following.
\begin{Proposition}
\label{prop6.4} In the above notations,
\begin{equation}
\label{6.24}
 H_0^* = (2 \cdot 4^m (m+1) (2m)!)^{-1}, \quad H_0 =
0.
\end{equation}
\end{Proposition}
With $H_0 = 0 $  we cannot find the asymptotic of $B^+(n) $  by
applying the same scheme which was successful in Sections 3-5. Notice
that in Section 5 we have the same difficulties in the case of
exceptional values of the coefficients of  $v. $  There $R^+ $  and
$R^- $ are  analogs of $H_0 $ (see (\ref{5.71}), (\ref{5.72}) and
(\ref{5.72a})). More precisely, if $v$ is given by (\ref{5.81}) then
$$ H^+_0 = \sum_{\xi \in X^+} h(\xi,0 ) = 0,  \quad
\text{if}   \; \begin{cases}    n \;  \text{is even and }  \;
 \cos \frac{\pi}{2} \sigma = 0,\\  n \; \text{is odd and} \;
 \sin \frac{\pi}{2} \sigma = 0,  \end{cases} $$
where $ \sigma^2 = -b^2/(4B), $  and
$$ H^-_0 = \sum_{\xi \in Y^-} h(\xi,0 ) = 0,  \quad
\text{if}   \; \begin{cases}    n \;  \text{is even and }  \;
 \cos \frac{\pi}{2} \tau = 0,\\  n \; \text{is odd and} \;
 \sin \frac{\pi}{2} \tau = 0,  \end{cases} $$
where $ \tau^2 = -a^2/(4A). $

 \vspace{3mm}

2. In this paper we consider only operators on the  interval
$[0,\pi].$ But let us mention that F. Gesztesy and V. Tkachenko
results \cite[Remark 8.10]{GT09} together with our examples from
Sections 3--5  show that {\em the Schr\"odinger operators
$$
Ly = - y^{\prime \prime}  + v(x) y, \quad x \in \mathbb{R},
$$
with potentials
$$
(1)  \quad v= ae^{-2ix} + be^{2ix}, \quad a,b \neq 0, |a| \neq |b|;
$$
$$
(2)  \quad v= ae^{-2ix} + be^{4ix}, \quad a,b \neq 0 ;
$$
$$
(3)  \quad v= ae^{-2ix} + Ae^{-4ix} + be^{2ix} + Be^{4ix}, \quad a,
b, A, B \neq 0, |A| \neq |B|,
$$
are not spectral operators of scalar type.}

\end{document}